\documentclass[letter,11pt]{article}

\usepackage{comment}
\usepackage[linktocpage=true]{hyperref}
\usepackage[margin=1.2 in, top=1 in, bottom= 1.2 in]{geometry}

\usepackage[english]{babel}
\usepackage{amsfonts}
\usepackage{mathrsfs}
\usepackage[mathscr]{euscript}
\usepackage{bbm}
\usepackage{latexsym}
\usepackage{math dots}
\usepackage{amssymb}
\usepackage{mathtools}

\usepackage{tikz}

\usepgflibrary{shapes.geometric}
\usetikzlibrary{calc}

\usepackage{enumitem}
\setlist{nolistsep}
\usepackage{amsthm}
\usepackage{relsize}
\usepackage{nicefrac}
\usepackage[capitalize]{cleveref}

\newtheoremstyle{plain}{3mm}{3mm}{\slshape}{}{\bfseries}{.}{.5em}{}
\newtheoremstyle{definition}{2mm}{2mm}{}{}{\bfseries}{.}{.5em}{}
\theoremstyle{plain}
	
\newtheorem{theorem}{Theorem}

\newtheorem{lemma}[theorem]{Lemma}

\newtheorem{corollary}[theorem]{Corollary}

\theoremstyle{definition}

\theoremstyle{plain}
\newtheorem*{namedthm}{\namedthmname}
\newcounter{namedthm}
	

\newcommand{\N}{\mathbb{N}}
\newcommand{\Z}{\mathbb{Z}}
\newcommand{\R}{\mathbb{R}}

\newcommand{\F}{\mathbb{F}}

\newcommand{\eps}{\epsilon}
\newcommand{\mdim}{\overline{D}}

\title{A discrete Marstrand type slicing theorem with the mass dimension.}
\author{Aritro Pathak}
\date{}
\begin{document}

\maketitle

\begin{abstract}
    \ \ \ \ Discretized versions of some central questions in geometric measure theory have attracted recent attention; here we prove a Marstrand type slicing theorem for the subsets of the integer square lattice. This problem is the dual of the corresponding projection theorem, which was considered by Glasscock, and Lima and Moreira, with the mass and counting dimensions applied to subsets of $\Z^{d}$. In this paper, more generally we deal with a subset of the plane that is $1$- separated, and the result for subsets of the integer lattice follows as a special case. We show that the natural slicing question in this setting is true with the mass dimension.
\end{abstract}
\text{Keywords: Marstrand slicing, Discrete slicing, mass dimension, counting dimension}
\section{Introduction and statement of results.}

\subsection{Dimensions of subsets of $\R^{2}$ and $\Z^2$, and the slicing theorem.}

\ \ \ \ \   Dimensions of fractal subsets of $\R^2$ are standard objects quantifying how ``large'' such fractals are \footnote{See for example, Chapter 4 of \cite{Mattila}}. Here we briefly review the basic notions, the classical Marstrand Slicing theorem, and then the corresponding natural analog of the slicing theorem for infinite $1$-separated subsets of $\R^2$.












Marstrand's slicing theorem is a general result in geometric measure theory, which states that for any Borel susbet $E\subset \R^2$,

\begin{theorem}\label{thm:slice1}
    $\text{dim}_H(E\cap l)\leq \text{max}(\text{dim}_H(E)-1,0)$ for \textit{almost all} straight lines $l$ in the plane.
\end{theorem}

Here, $\text{dim}_H(F)$ denotes the Hausdorff dimension of any set $F\subset \R^{2}$.

 Moreover, the bound on the right hand side of \cref{thm:slice1} is the smallest possible for the statement to hold in general; i.e. for any small $\eps>0$, it is possible to construct a Borel subset $E_\eps \subset \R^2$ such the set of lines satisfying $\text{dim}_{H}(E_\eps \cap l) > \text{max(dim}_H(E)-1-\eps,0)$ has positive Lebesgue measure. Here the set of lines in $\R^2$ is a two parameter set, where the parameters can be chosen to be the slope and the $y$ intercept. 

In fact, \cref{thm:slice1} is a consequence of the following stronger statement, which is stated in Theorem 1.6.1 in \cite{Bishop}.

\begin{theorem}\label{thm:slice2}
Let $E\subset \R^2$ and let $E_x=\{ y:(x,y)\in E\}$. If $\text{dim}_H(E) \geq 1$, then $\text{dim}_H(E_x)\leq \text{dim}_H(E)-1$ for Lebesgue almost every $x$, and if $\text{dim}_H(E)<1$, then the slices $E_{x}$ are empty for Lebesgue almost every $x$.
\end{theorem}

Clearly, \cref{thm:slice2} implies \cref{thm:slice1}: we can simply reorient the axes in the plane and apply \cref{thm:slice2} to all the slices perpendicular to the rotated $x$ axis, and in the process get the statement for all possible lines in the plane. \cref{thm:slice2} follows from basic considerations of the definition of the Hausdorff dimension.

\bigskip

In this paper, we are concerned with $1$-separated sets in $\R^2$\footnote{These are sets whose every pair of points has distance at least $1$}, and the natural notion of dimension is in analogy with the box dimension of arbitrary subsets of $\R^{2}$. Dimensions of subsets of $\Z$ or $\Z^{d}$ have been studied by Barlow and Taylor \cite{Barlow1, Barlow2}, Naudts \cite{Naudts1, Naudts2}, Iosevich, Rudnev, and Uriarte-Tuero \cite{Iosevich}, Lima and Moreira \cite{Moreira}, and Glasscock \cite{DG}. In these papers parallels are drawn between this infinite discrete setting and the classical setting of arbitrary subsets of $\R^{d}$. We define the mass and counting dimensions of $1$-separated subsets of $\R^{2}$ below. The definitions generalize in the obvious way to $1$-separated subsets of $\R^{d}$.

The mass and counting dimensions of any $1$-separated set $E\subset \R^2$ are respectively defined as: 

\begin{equation}\label{eq:dimen}
    \mdim(E)=\limsup\limits_{l\to \infty}\frac{\log|E\cap [-l,l]^{2}|}{\log(2l)}, \ D(E)=\limsup\limits_{||C||\to \infty}\frac{\log|E\cap C|}{\log||C||}.
\end{equation}

Here, $|\cdot |$ represents the cardinality of a set. The limit supremum for the counting dimension is taken over all possible cubes, with the side length $||C||$ of the cubes going to infinity. Thus in considering the mass dimension, we always consider cubes with centers at the origin, while for the counting dimension, the limit supremum is taken over arbitrary cubes whose lengths are going to infinity\footnote{Even in $\R^{2}$ we loosely refer to $C$ as a cube.}. Thus in particular, we always have $\mdim(E) \leq D(E)$.

Two prototypical examples of ``fractal" sets in the integers are values of polynomials with integer coefficients at integer arguments, and restricted digit Cantor sets. The latter are integer Cantor sets with restrictions on the digits that are permissible to be used. These examples have been studied in detail in \cite{DG},\cite{Moreira}, with the mass and counting dimensions as the natural dimensional parameters.

In our case, in place of lines it is natural to consider tubes of width 1, that extend to infinity. For simplifying our arguments, we restrict our set of points to lie only in the first quadrant. Qualitatively nothing changes, the set still extends to infinity inside the first quadrant, and we modify the definition of the mass dimension to:

\begin{equation}
    \mdim(E)=\limsup_{l \to \infty} \frac{\log|E \cap [0,l]^{2}|}{\log(l)}
\end{equation}

It is clear that this definition is equivalent to the earlier one; we only got rid of an additive term of $\log 2$ from the denominator, and since $l\to \infty$ the definitions are equivalent. 

For $E\subset \N$, the definition reduces to 

\begin{equation}
    \overline{D}(E)= \limsup\limits_{N\to \infty} \frac{\log|(E\cap \{1,2,...,N\})|}{\log(N)}
\end{equation}

In place of $1$-separated sets, one can consider an arbitrary $\delta$-separated set, and then consider tubes of width $\delta$ in our statements, and the arguments would essentially remain the same. Our motivation for considering a general $1$-separated set is to apply our results to the integer grid which is naturally $1$-separated. Finite $\delta$-separated sets of arbitrary dimensions have been considered by Katz and Tao \cite{KT} when they define $\delta$- discretized analogs of the Falconer distance conjecture, the Erd\H{o}s Ring conjecture and the Furstenberg set. Analogs of Marstrand's classical projection theorem for a finite $\delta$-discretized setting have been obtained by Shmerkin \cite{Shmerkin2,Shmerkin3} and Rams \cite{Rams}. This essentially involves using the Riesz energy and a Tchebysheff type argument. The projection theorems in our general setting, using both the mass and counting dimensions, have been earlier proven by Lima and Moreira \cite{Moreira}, and Glasscock \cite{DG}. \footnote{We could in principle consider an arbitrary point set in $\R^{2}$, but when we have some limit points in our set, then the dimension of the set goes to infinity. For $1$-separated sets in the plane, the dimension is at most 2.} 

In this paper, we prove a Marstrand type slicing theorem for the mass dimension. 
With an elementary Tchebysheff inequality argument, we first show that our result holds true in a weaker ``asymptotic'' sense, a notion we define in the next section. This weaker result only depends on the arithmetic properties of our sets in question, and we find an elegant analogous result in the finite field setting as well, which we state in the next subsection. Working in a finite field makes clear how one then tackles the corresponding problem in $\R^2$.

The classical Marstrand projection and slicing theorems do not give any information about the dimensions of specific projections or specific slices. In some situations, for certain types of sets it is known that there is no dimension drop in a specific projection, \footnote{see for example \cite{Shmerkin2,Kauffman}}. It was a longstanding conjecture of Furstenberg that for $p,q$ multiplicatively independent, (i.e $\frac{\log p}{\log q}$ being irrational), if we consider two subsets $A,B$ of $[0,1)$ invariant under the $\times p$ and $\times q$ maps respectively, then every slice of $A\times B$ satisfies the statement of the Slicing theorem.

Equivalently, if $g$ is any affine transformation, under the assumptions stated above, we have:

\begin{equation}
    \text{dim}_{H}(A\cap g(B))\leq \text{max}(\text{dim}_{H}(A)+\text{dim}_H(B)-1,0).
\end{equation}

Recently there have been two proofs of this celebrated conjecture $\cite{Shmerkin1,Wu} $. In the integer setting there are longstanding analogous transversality questions for the particular case of restricted digit Cantor sets. A discussion on this is found in the introduction in \cite{GMR}; there is a known conjecture \cite{OEIS} about the set of integers $A$ that can be written in bases $2,3,4,5$ with digits only $0,1$ being $\{0,1,82000\}$. This is obviously a question about the intersections of multiple restricted-digit-Cantor-sets, which are sets respectively invariant under $\times 2,\times 3,\times 4,\times 5$ maps. Recently, Burrell and Yu \cite{Yu} showed that the above defined set $A$ is such that $|A \cap [0,N]|\leq C_\eps N^{\eps}$ for any positive $\eps>0$, which shows that the intersection set has mass dimension 0. This is another illustration where the mass dimension is the natural dimensional quantifier. A well known conjecture of Erd\H{o}s \cite{Erdos} which states that except for finitely many natural numbers $n$, the integer $2^{n}$ contains $1$ in its base $3$ expansion, can also be naturally formulated in terms of intersections of specific digit Cantor sets. 

While this illustrates longstanding interest in upper bounds to the dimensions of intersections of specific structured subsets of the integers, or more generally the intersection of one such subset and with the other under a natural affine map, a more fundamental question is whether such a natural upper bound can be obtained for almost all such affine intersections in a Lebesgue sense, and this is the Marstrand type slicing theorem that we prove in this paper.

We also remark that there is no analogous version to \cref{thm:slice2} stated earlier, for our slicing theorem for the mass dimension. In fact, we will easily construct examples of a set $E$ where for a given projecting direction, every tube perpendicular to this projecting direction has mass dimension greater than $\mdim(E)-1$. However, \cref{thm:slice1} holds true, due to the fact we prove later, that if we consider a ray of tubes centered at some point, Lebesgue almost all these tubes will have mass dimension less than or equal to $(\mdim(E)-1)$.

For an analogous version of our problem with finitely many points, the Szemeredi Trotter theorem gives an upper bound to the number of incidences between points and lines. This was generalized recently in \cite{Guth} to the setting of $\delta$-tubes and $\delta$-balls, with some restrictions on the spacing of the tubes. In our case, we are dealing with an infinitary asymptotic analog of this problem with $\delta=1$.

\bigskip 

\subsection{Statement of results}

We first prove a result in finite fields, where the result holds for almost every line of the finite field in an asymptotic weaker sense. We replicate the same argument in the real plane, and get a corresponding result in the real plane in the weak asymptotic sense. In the real plane, care needs to be taken while dealing with a double limit when taking the length of the grid of points as well as the parameter set of the lines to both go to infinity, which is a subtlety averted in the finite field problem due to the fact that the set of lines as well as the set of points both grow as $p^{2}$. These results are simpler but interesting in that they are obtained by purely arithmetic manipulations using Fubini and Tchebycheff type arguments; however our strongest result, \cref{thm:strongerslicing} depends on the topology of the plane.

Consider the finite field $\mathbb{F}_{p}$ with $p$ prime; for subsets $B \subset \mathbb{F}_{p}^{2}$, we modify the definition of the mass dimension in the natural way:

\[\overline{D_p}(B)=\frac{\log|B \cap \{1,2,\ldots,p\}^{2}|}{\log p} \] 

For any $(u,v) \in \mathbb{F}_{p}^{2}$, we have the line $\ell_{(u,v)}:=\{(x,y)\in \mathbb{F}_{p}^{2}|y=ux+v\}$. In the real plane the analog of this line is $l_{(u,v)}:=\{(x,y)\in \mathbb{R}^{2} |y=ux+v\}$ with $u,v\in \R $.\footnote{We also define $\lfloor l_{(u,v)}\rfloor :=\{(x,y)\in \mathbb{R}^{2} | y=\lfloor ux+v \rfloor \}$} For any subset $E\subset \Z^2 $, we define in the real plane the set $E \cap \lfloor l_{(u,v)}\rfloor := \{ y\in \R| (x,y)\in E, y=\lfloor ux+v \rfloor \}$. In particular, if there are multiple points of $E$ in any horizontal strip of the set $\{(x,y)\in \R^2 |y=\lfloor ux+v \rfloor\}$, we only have a count of 1 in $(E\cap \lfloor l_{(u,v)} \rfloor)$, for all of these points.

\begin{theorem}\label{thm:weakfinite }
 \emph{For all $E \subseteq \mathbb{F}_{p}^2$, the set $U$ of parameters $(u,v) \in \mathbb{F}_{p}^{2}$, so that}
\[\overline{D_p} \big(E \cap  \ell_{(u,v)} \big) \leq \max \big(0, \overline{D}(E) - 1 \big)+o_p(1)\] \emph{is such that $\lim\limits_{p\to \infty} \frac{|U \cap \mathbb{F}_{p}^{2}|}{p^{2}}=1$,  where the limit is taken along the sequence of primes.}

\end{theorem}

This gives the following corollary:

\begin{corollary}\label{corollary:weakfinite}
\emph{For all $A, B \subset \mathbb{F}_{p}^2$, the set $U$ of parameters $(u,v)\in \mathbb{F}_{p}^2$, so that }
\[\overline{D_p} \big(A \cap  (u B + v)  \big) \leq \max \big(0, \overline{D_p}(A) + \overline{D_p}(B) - 1 \big)+o_p(1)\] \emph{is such that $ \lim\limits_{p\to \infty} \frac{|U \cap \mathbb{F}_{p}^{2}|}{p^{2}}=1$, where the limit is taken along the sequence of primes.}
\end{corollary}

We now state the corresponding results in the real plane.

\begin{theorem}\label{thm:weakreal}
\emph{For all $E \subseteq \N^2$, the set $U$ of parameters $(u,v)$, with $u,v>0$, so that}
\[\overline{D} \big(E \cap \lfloor l_{(u,v)}\rfloor \big) \leq \max \big(0, \overline{D}(E) - 1 \big)\] \emph{is such that $\lim\limits_{M\to \infty} \frac{|U \cap [0,M]^{2}|}{M^{2}}=1$. }

\end{theorem}

This gives the following corollary for the Cartesian grid:

\begin{corollary}\label{corollary:weakreal}
 \emph{For all $A, B \subseteq \N$, the set $U$ of parameters $(u,v)$, with $u,v>0$, so that }
\[\overline{D} \big(A \cap \lfloor u B + v \rfloor \big) \leq \max \big(0, \overline{D}(A) + \overline{D}(B) - 1 \big)\] \emph{is such that $ \lim\limits_{M\to \infty} \frac{|U \cap [0,M]^{2}|}{M^{2}}=1$.}
\end{corollary}

We prove these results in Section 2.\footnote{These conclusions above have the same structure as the hypothesis in the geometric Ramsey theory results in \cite{Katznelson} and \cite{Bourgain} where sets with such asymptotic full density in the plane are shown to always contain certain geometric configurations.}

We now state the main result of this paper, which we will prove in Section 5. Here for convenience, the parametrization $(u,v)$ changes from that in the earlier theorems.

The tube $t_{u,v}$ is explicitly described as
\[t_{u,v} = \left\{(x,y) \in \R^2 \ \middle| \ -\frac 1u x + v \sqrt{ 1 + \frac 1{u^2}} < y \leq -\frac 1u x + (v+1) \sqrt{ 1 + \frac 1{u^2}} \right\}.\]

It is easily seen that this is a tube whose width is exactly of length 1. The perpendicular line to this tube has slope $u$. Consider the line perpendicular to this tube passing through the origin. The coordinate $v$ gives us the displacement of the tube from the origin, along this perpendicular line.

\begin{theorem}\label{thm:slicing}
Let $E \subseteq  \R^2$ be a $1$-separated set of mass dimension $\mdim(E)$. Then in the Lebesgue sense, for almost every tube $t_{u,v}$ of width $1$, slope $u$, and displacement $v$ along the projecting line, we have that $\mdim(E \cap t_{u,v}) \leq \text{max} (0, \mdim(E)-1)$.
\end{theorem}

We choose the tube in this manner, closing the upper edge and keeping the lower edge open, since we eventually wish to apply this theorem to the broken line of the form $\lfloor l_{u,v}\rfloor=\{(x,y): y=\lfloor ux+v \rfloor\}$ with $u>0$, and so for such a tube (whose width is less than 1 but the vertical intercepts are 1), we should include the top edge and keep the bottom edge open.

We will prove a stronger statement than \cref{thm:slicing}:

\begin{theorem}\label{thm:strongerslicing}
Let $E \subseteq  \R^2$ be a $1$-separated set of mass dimension $\mdim(E)$. Then for all $v \in \R$, for Lebesgue-a.e. $u \in \R_+$, 
\[\mdim(E \cap t_{u,v}) \leq \mathrm{max} (0, \mdim(E)-1).\]
\end{theorem}

Upon integrating over all $v\in \R$, \cref{thm:strongerslicing} implies \cref{thm:slicing}.

We will get the following corollary specific to the Cartesian grid:

\begin{corollary}\label{corollary:finalcor}
Let $A, B \subset  \N$. For almost every $u, v \in \R^2$, 
\[\mdim \big(A \cap \lfloor u B + v \rfloor \big) \leq \max \big(0,\mdim(A)+\mdim(B)-1 \big).\]
\end{corollary}

The formulation in the above corollary is how we immediately apply the theorem to specific subsets of $\N$ such as polynomials or restricted-digit-Cantor-sets.

In essence, the purpose of dealing with $1$-separated sets and tubes is to be able to make a statement about the dimension of the lines within the integer grid, as above. The whole argument also works with $\delta$-separated sets and $\delta$ tubes.

While we prove \cref{thm:strongerslicing} in $\R^2$, the analogous result in higher dimension should follow in the same way by contradiction, and our method of integration should apply in essentially the same manner. We also do not study the question of how large the dimension of the set of exceptional tubes can be, for any given set. In Examples 1 and 2 of Section 3, we construct sets where this set of exceptional tubes is of dimension one, but it remains open to construct examples where the dimension of this set of exceptional tubes is greater than 1.

Since the preparation of this manuscript, there has been further work \cite{Pathak} showing that similar Marstrand type slicing results are not true when one considers the counting dimension (as defined in \cref{eq:dimen}) in place of the mass dimension that has been used throughout this paper. In fact, it is shown that not even the weaker \cref{thm:weakreal} holds true when one considers the counting dimension and there are counterexamples.

\bigskip

\subsection{Outline of proof of \cref{thm:strongerslicing}}

We now outline the proof of the main result of this paper, \cref{thm:strongerslicing}. The proof runs by contradiction. Assume to the contrary that there is some set $E$, some $v_0 \in \R$ such that there is a positive Lebesgue measure set of $u$'s such that $\mdim(E \cap t_{u,v_0}) > \text{max} (0, \mdim(E)-1)$. This clearly implies that there is some $\epsilon>0$ such that there is some positive measure  set $U$ of $u$'s such that $\mdim(E \cap t_{u,v_0}) > \text{max} (0, \mdim(E)-1)+\eps$.

The idea is to restrict $U$ to some small interval, say $U'$ and without loss of generality consider $v_0=0$. Thus all the tubes parametrized by values in $U'$ pass through the origin, and lie within a small cone, corresponding to the interval. The area of a cone grows as the square of the length of the edge of the cone, and it is thus a two dimensional object (unlike a vertical strip of finite horizontal width, which is essentially a one dimensional object). A positive Lebesgue measure of the tubes within this cone have an exceptionally high dimension. So the idea is to carefully `integrate' the number of points across each of these tubes, so that the cone in itself has an exceptionally high number of points and thus has dimension greater than $E$ itself, giving a contradiction. 

In order to perform this `integration' within the cone, we will show that for each of the exceptional tubes, there exist heights $H_{l,t}$, indexed by the integers, so that in the upper half of the tube till these heights in the range $[H_{l}/2, H_{l}]$, we have an exceptional number of points. We then carefully choose a positive Lebesgue measure subset of the tubes so that these special heights of all these tubes in this specific subset are in a quantifiably similar height range. Now within this subset of tubes, we choose a well spaced set of tubes so that the upper halves of all these tubes are mutually disjoint, so there is no  double counting of any point. We add up all the points in the upper halves of all these well spaced tubes of this subset, which are all this the approximately similar height range. This would imply that the intersection of the cone with a suitable square whose one side extends vertically from the $x$ axis till anywhere in this approximate height range, in itself has an exceptionally large number of points. This process can be repeated for a countably infinite number of such height ranges. This then shows that the intersection of the cone with $E$ in itself has dimension greater than $\mathrm{dim}(E)$ which gives us the contradiction.

To simplify the argument, without loss of generality, we choose the cone to be pointed vertically upward, when we prove the result in Section 4. This same argument does not work if instead of the cone, we work with an exceptional set of tubes that belong within a fixed strip, i.e. with fixed $u$ coordinates but varying $v$ coordinates. The reason essentially boils down to the fact that the strip under consideration is essentially a one dimensional object, and it cannot guarantee that we find a subset of $E$ with exceptionally large dimension. In fact in Section 3, we give a simple example of a set where every single tube with a fixed $u$ coordinate but varying $v$ coordinates has an exceptionally high dimension.

\bigskip

In Section 2 of the paper, we prove the asymptotic results in finite fields and in the real plane. In Section 3, we cite several different sets to illustrate that \cref{thm:strongerslicing} is optimal, and that the inequality cannot be made stronger.\footnote{Note that this does not prove that the inequality in the weaker result that is  \cref{thm:slicing} is optimal} In Section 4, we prove \cref{thm:strongerslicing} and then \cref{corollary:finalcor}. Throughout the paper, as used earlier, we use $D_p$, $\overline{D}$ to denote respectively the mass dimension in the finite filed $\F_{p}^{2}$, the mass dimension in the real plane $\R^{2}$.

\section{Asymptotic results for $\mathbb{F}_{p}^2$ and $\Z^2$}

First we prove \cref{thm:weakfinite } for the case of finite fields, and then \cref{corollary:weakfinite} .  

\begin{proof}[Proof of \cref{thm:weakfinite }] Consider the function $k(p)=\log p$. We show with the basic Fubini and Tchebysheff type arguments that for at least $p^{2}(1-\frac{1}{k(p)})$  of all the possible lines $l_{u,v}$ in $\mathbb{F}_{p}^{2}$ we have 
$|E\cap l_{u,v}|\leq \frac{k(p)}{p}|E|$. Thus, from here we conclude our result for the fraction $(1-\frac{1}{\log p})\big( \rightarrow 1 \ \text{as} \ p\to \infty \big)$ of all possible lines in $\F_{p}^{2}$. This is clearly a satisfactory "almost every" line description as $p \to \infty$.

Given any pair $(x,y)\in E$, clearly for any $u\in \mathbb{F}_{p}$, there is exactly one $v\in \mathbb{F}_{p}$ so that $y=ux+v$, and so for each pair $(x,y)\in E$, we have exactly $p$ possible pairs so that $y=ux+v$ (and furthermore, we cannot have two pairs of the form $(u,v_{1}), (u,v_{2})$ with $v_{1}\neq v_{2}$ in this set of pairs, nor two pairs of the form $(u_{1},v),(u_{2},v)$ with $u_{1}\neq u_{2}$, as $p$ is a prime). We have $|E|$ possible points $(x,y)\in E$, and for each of them we have $p$ possible pairs $(u,v)$ with the above property. Since there are $p^{2}$ possible pairs $(u,v)\in \mathbb{F}_{p}^{2}$, this implies that on average over $\mathbb{F}_{p}^{2}$, the set $|E\cap l_{u,v}|$ has $\frac{|E|p}{p^2}=\frac{|E|}{p}$ number of elements, with $\sum\limits_{(u,v)\in \mathbb{F}_{p}^{2}}|E\cap l_{u,v}|=|E|p$.



We see that at least $p^{2}(1-\frac{1}{k(p)})$ of the pairs are such that $|E\cap l_{u,v}|\leq \frac{k(p)|E|}{p}$. If not, at least  $\frac{p^{2}}{k(p)}$ of the pairs are such that $|E\cap l_{u,v}|\geq \frac{k(p)|E|}{p}$ and these alone would sum to $\geq |E|p$, and we have a contradiction. 

It's easy to see upon taking logarithms, that the above means that for a $(1- 1/k(p))$ proportion of the pairs $(u,v)$,
\[\overline{D_p} \big(E \cap l_{u,v})\big) \leq D_p(E) - 1 + \frac{\log k(p)}{\log p}.\]
Now since $k(p)=\log p$, as alluded to in the beginning of the proof, this shows that for almost all pairs $(u,v)$ in the asymptotic sense, $D_p \big(E\cap l_{u,v})\big) \leq max(0,D_p(E)-1) + o_{p}(1)$. 
\end{proof}

Note that the set $(A \cap (uB+v))$ is the same set  as the intersection of the affine line $\ell_{u,v}$ which is $\{(x,y)\in \mathbb{F}_{p}^{2}: y=ux+v \}$, with $A \times B$. The proof of \cref{corollary:weakfinite} thus follows by simply considering in place of $E$ the Cartesian grid $A\times B$.
\bigskip

Now we prove \cref{thm:weakreal} and then \cref{corollary:weakreal}.

\begin{proof}[Proof of \cref{thm:weakreal}]
Suppose $E\subseteq \N^{2}$. For $N \in \N$, let $E_N = E \cap \{1, \ldots, N\}^{2}$,  let for $(a,b)\in E$,
\[T_{(a,b)} = \{(u,v) \in [0,M]^2 \ | \ a = \lfloor ub + v \rfloor \}.\]

It is easily seen that when $u> 1$, each horizontal strip of the broken line segment can contain at most one point belonging to $E_{N}\subset \mathbb{N}^{2}$, and so we have the basic counting identity:
\[\big| E_{N} \cap \lfloor l_{(u,v)} \rfloor \big| = \sum_{\substack{(a,b) \in E_N}} \chi_{T_{(a,b)}}(u,v)\]

where $\chi_{T_{(a,b)}}$ is the indicator function of the set $T_{(a,b)}$.

Otherwise when $0<u\leq 1$, we clearly have:

\[\big| E_{N} \cap \lfloor l_{(u,v)} \rfloor \big| \leq \sum_{\substack{(a,b) \in E_N}} \chi_{T_{(a,b)}}(u,v)\]

More generally we can simply use this inequality above for all $u> 0$.

We employ the double counting plus Tchebysheff inequality technique in this problem, similar to the finite field case.

Consider for some large $M>0$, the integral,

\[\frac 1{M^2} \iint_{(0,M]^{2}} \sum_{\substack{(a,b) \in E_N}} \chi_{T_{(a,b)}}(u,v) \ dudv =  \frac 1{M^2} \sum_{\substack{(a,b) \in E_N}} \iint_{(0,M]^{2}}  \chi_{T_{(a,b)}}(u,v) \ dudv    \]

Consider some fixed $(a,b)\in E_N$. If the integral on the right side, for this given $(a,b)\in E_N$, were to be taken over the entire $(u,v)$ plane, and not just restricted to $(0,M]^2$, then we observe that for any $u_0\in \mathbb{R}$,  $\chi_{T_{(a,b)}}(u_0,v)=1$ whenever $a=\lfloor u_0 b+v \rfloor$, i.e. whenever $a\leq u_0 b +v <a+1$, or in other words whenever $v$ lies in the interval $[a-u_0 b,a+1-u_0 b)$ of width 1. 

When $v$ is bounded in the interval $(0,M]$ as is the case here, then given any $u_0 \in [0,M]$, clearly the set of  values of $v$ so that $\chi_{T_{(a,b)}}(u_0,v)=1$, has width less than or equal to 1.  
Thus on the integral on the right, we just integrate over the $u$ variable, noting that for such a fixed value of $u$, the integrand takes the value $1$ in a set of width at most one.

Thus the upper bound comes to,

\[ \frac{1}{M^{2}}\iint_{(0,M]^{2}} |E_N \cap \lfloor l_{(u,v)} \rfloor|du dv\leq  \frac 1{M^2} \sum_{\substack{(a,b) \in E_N}} \iint_{(0,M]^{2}}  \chi_{T_{(a,b)}}(u,v) \ du dv \leq  \frac{1}{M} |E_{N}|. \] 

The Tchebysheff argument now applies as in the finite field case, and clearly the area of $(0,M]^{2}$ which consists of pairs $(u,v)$ such that $|E_{N} \cap \lfloor l_{{(u,v)}})\rfloor|\leq \frac{k(M,N)}{M}|E_N|  $ is at least $M^{2}(1-\frac{1}{k(M,N)})$, where $k(M,N)$ is any positive function of $M,N$ which goes to infinity sublinearly in both $M,N$.

Now we would take the limit of $N,M$ both going to infinity, seemingly having freedom over both $M,N$ while doing so. However, if we always choose $M\gg N$ in the process of taking the simultaneous limits $M,N \to \infty$, then in each step of the limiting process, except a small fraction $O(\frac{N}{M})$, all of the other lines would have $v\gg N$. In that case each of these lines would have no intersection with the set $E_N$ in each step of the limiting process, and this is an incorrect way of taking the limit\footnote{if $M\gg N^{2}$ we would conclude that for the majority of the lines $\lfloor l_{u,v} \cap E \rfloor$ has zero cardinality, since the upper bound is $\frac{k(M,N)}{M}E_N\leq \frac{k(M,N)}{M}N^{2} \to 0$ in this case if $M\gg N^{2}$.}. On the other hand, while taking the limit if we always choose $N\gg M$, then the upper bound for $|E_N \cap \lfloor l_{(u,v)} \rfloor|$ becomes far weaker than optimal, although there is no counting error.

Optimally we take $M=\frac{N}{\alpha}$, where $\alpha>1$ and $\alpha=O(1)$. In that case, all the lines under consideration during each step of the limiting process intersects $E_N$, and when we take the limit of $N \to \infty$, we would get the correct limiting value for each pair $(u,v)$ with $u>1$ and $ v>1$. We consider $k(M,N)=\sqrt{\log M \log N}$. For the purpose of taking the limit, we might as well take $\alpha=1$, and the result we have is that for asymptotically almost every line $\lfloor l_{(u,v)} \rfloor$, 


\[ \frac{\log(|E_N \cap \lfloor l_{(u,v)} \rfloor|)}{\log N} \leq \frac{\log(\log N)}{\log (N)} +\frac{\log |E_N|}{\log N}-1  \]

and upon taking the limit as $N\to \infty$ we would get 

\[ \overline{D} \big(E \cap \lfloor l_{(u,v)} \rfloor \big) \leq \max \big(0, \overline{D}(E) - 1 \big)\].
\end{proof}

Now, \cref{corollary:weakreal} is a special case of \cref{thm:weakreal}, specializing to a rectangular grid $A\times B \subset \N^{2}$, since it is clear that $(A\times B) \cap \lfloor l_{(u,v)} \rfloor=A\cap \lfloor uB+v \rfloor $.

\section{Mass dimensions of specific sets and their slices.}

\begin{enumerate}

    \item We first give a trivial example illustrating that the statement of \cref{thm:strongerslicing} would not hold if the parameters $u,v$ were interchanged, i.e. one can construct a set $E$ so that  for a fixed $u$, for every $v>0$, the tubes $t_{u,v}$ have dimension greater than $\overline{D}(E)-1$. Simply consider the line $y=ux$ for any $u>0$, in the first quadrant, and put points at intervals of length 1 on this line. This is clearly a set of mass dimension one. Moreover, every tube  $t_{(-\frac{1}{u},v)}$ \ , where $v$ lies in the interval $[0,1)$ also intersects every point on this line, and thus also has mass dimension exactly 1. 

    \item 
We give a second less trivial example illustrating that  \cref{thm:strongerslicing} would not hold if the parameters $u,v$ were interchanged in the statement of the theorem.  Without loss of generality, consider $u=\infty$, that is, all the horizontal tubes of width 1. Consider the set $E:=\{(m^{2},n)| m\in \Z_+, 0< n\leq m \}$. Thus inside a square box with one side of length $N^{2}$ along the $x-$axis, there are $1+2+...+N=\frac{N(N+1)}{2}$ many points belonging to the set $E$. It is not hard to see that the mass dimension of the set is given by:
    
    \[\lim_{N\to \infty} \frac{\log(N(N+1)/2)}{\log(N^{2})}=1\]
    
    Each horizontal tube is initially empty, but then contains a point for every square value of the $x-$coordinate, for every sufficiently large square. This is the number of points in the intersection of the tube with the square $[0,N]^{2}$. Thus for sufficiently large $N$, there are effectively $\sqrt{N}$ many points in each tube, and thus each tube has mass dimension $1/2$, and thus each of these tubes is an exceptional tube.

     However, for a fixed value of $v_{0}$, we cannot have a positive Lebesgue measure set of values of $u$ so that the tubes $t_{(u,v_0)}$ are exceptional, as stated in \cref{thm:strongerslicing}. This is proven in Section 4.
     
\end{enumerate}     
    
    \bigskip
    
    Next we show several examples (Examples 3-6) where equality is attained in the statement of \cref{thm:strongerslicing}, illustrating that the inequality cannot be strengthened in general. Example 3, which we state below, illustrates that if one attempts to construct a 1 dimensional set and then find a ray of tubes whose $v$ coordinates are in a set of positive Lebesgue measure, and each tube with positive mass dimension, one will fail and the ray of tubes will all end up having 0 dimension. On the other hand, going from Example 4 to Example 6 is a prototypical illustration of the fact that if one attempts to construct a ray of tubes whose $u$ coordinates are in a set of positive Lebesgue measure, with each tube having an exceptionally high dimension, one will fail and that eventually the set itself will be forced to have a high enough dimension so that statement of \cref{thm:strongerslicing} holds true. We also note that for a tube of width 1 in any direction, the intersection of the integer lattice with a length $m$ of this tube has between $m$ and $2m$ many points. 
    
    \bigskip
    
\begin{enumerate}[resume]    
    
    \item The following example exhibits a set $E$ of mass dimension 1, a fixed value $v_0$, and an open set of $u$'s so that every slice by tubes $t_{u,v_0}$ has mass dimension exactly 0.
    
    \begin{figure}
\centering

\begin{tikzpicture}[x=1cm, y=1cm, semitransparent, scale=0.6]

\draw [->] (0,0) -- (0,6) node[above] {$y$} ;
\draw [->] (0,0) -- (6,0) node[right] {$x$} ;
\draw [thick, dotted] (0,0) -- (5.5,5) ;
\draw [thick, dotted] (0,0) -- (5,6)  ;
\draw [thick,dashed] (1,1) -- (1,1.2) ;
\draw [thick,dashed] (1,1.2) -- (1.35,1.2) ;
\draw [thick,dashed] (1.35,1.2) -- (1.35,1.6) ;
\draw [thick,dashed] (1.35,1.6) -- (1.7,1.6) ;
\draw [thick,dashed] (1.7,1.6) -- (1.7,2) ;
\draw [thick,dashed] (1.7,2) -- (2.2,2) ;
\draw [thick,dashed] (2.2,2) -- (2.2,2.7) ;
\draw [thick,dashed] (2.2,2.7) -- (3,2.7) ;
\draw [thick,dashed] (3,2.7) -- (3,3.6) ;
\draw [thick,dashed] (3,3.6) -- (4,3.6) ;
\draw [thick,dashed] (4,3.6) -- (4,4.8) ;
\draw [thick,dashed] (4,4.8) -- (5.3,4.8) ;
\draw [thick,dashed] (5.3,4.8) -- (5.3,5.4) ;
\end{tikzpicture}
\caption{The set used in Example 3 depicted in the figure with the dark dashed line. It has alternate horizontal and vertical parts within the annual region depicted in the figure, all points of the set belonging to the integer lattice}
\end{figure}
    
    Consider the slanted closed cone with vertex at the origin, the right edge being the line $y=x$ and the left edge being the line $y=x\tan(\pi/4 +\delta)$ for some small $\delta$. Begin at the point $(1,1)$ and consider all the integer points on the vertical line starting from $(x_1,y_1)=(1,1)$ to $(1,\tan(\pi/4 +\delta))$. From here we take a `right turn' along the horizontal line and consider all the integer points on it, till it hits the line $y=x$, and then again consider the alternate vertical and horizontal lines (with a succession of `up' and `right' turns). This `zig-zag' growing sequence of points is our set $E$. The growth of this set is linear, and it is not hard to verify that the mass dimension of this set is exactly $1$. Consider the points $(x_{n},y_{n}), n\in \mathbb{N} $, where the horizontal segments of our `zig-zag' line intersects the line $y=x$ ( and thus always $x_{n}=y_{n}$).
    
    Consider all the tubes of width 1, whose right edge passes through the origin and with these right edges being within the cone. It is clear that for each of these tubes, for any given $n\in \mathbb{N}$ there are two points of the set $E$ that lie within this tube, corresponding to the horizontal and vertical lines that emanate from the point $(x_n,y_n)$. So we need to check how fast the coordinates of the points $(x_n,y_n)$ grow.
    
    With some elementary trigonometry, it can be verified that the growth can be written iteratively, for some small positive constants $A= (\cot(\pi/4)-\cot(\pi/4+\delta)) >0,B=(\tan(\pi/4 +\delta) -\tan(\pi/4))>0$ independent of $n$, in the following way:\footnote{From the point $(x_n,y_n)$, till we hit the line $y=x\tan(\pi/4 +\delta)$ the $x$-coordinate remains constant while the increment in the $y$-coordinate is exactly $x_n\cdot (\tan(\pi/4 +\delta)-\tan(\pi/4)))$. Further from this point, the next time we hit the line $y=x$ the $y$-coordinate remains constant while the $x$-coordinate increased by exactly $y_n\cdot(\cot(\pi/4)-\cot(\pi/4+\delta)$. }

\begin{center}    
    $\begin{pmatrix} x_{n+1} \\ y_{n+1} \end{pmatrix} =\begin{bmatrix} 1&A\\ 0&1 \end{bmatrix} \begin{bmatrix}1&0\\B&1 \end{bmatrix} \begin{pmatrix} x_n \\ y_n\end{pmatrix}= \begin{bmatrix}1+AB & A \\ B &1\end{bmatrix} \begin{pmatrix}x_n \\ y_n\end{pmatrix}$
\end{center}    

\newpage   
   Thus, 
   
\begin{center}
   $ \begin{pmatrix} x_{n}\\y_{n} \end{pmatrix}= \begin{bmatrix}1+AB & A\\ B & 1 \end{bmatrix}^{n} \begin{pmatrix}1\\ 1\end{pmatrix}$
\end{center}

    The eigenvalues of the matrix above are $\lambda_1=1+ \frac{1}{2}(\sqrt{AB(AB+4)}+AB)>1$, $\lambda_2=1-\frac{1}{2}(\sqrt{AB(AB+4)}-AB)<1$, and thus we have hyperbolicity, and conclude that $x_{n}$ has an additive term with a factor that grows like $\lambda_{1}^{n} $, that is the growth is exponential. Thus till the height $x_{n}$, for each of the tubes within the cone, we have $\sim 4O(\log(x_{n}))$ elements within the tube. The additional factor of 2 comes from the fact that each tube within the cone contains two points `near' each $x_n$, one each from the horizontal and vertical strips which intersect the point $(x_n,y_n)$.   Thus the mass dimension of each tube is exactly:
    
    \[ \lim_{n\to \infty} \frac{\log(2\log(x_{n}))}{\log(x_{n})}=0 \]
    
    \item This is a simple example of a set $E$ that has mass dimension $3/2$ while a set of tubes of the form $t_{(u,v)}$ with fixed $u$, and the $v$ parameter lying in an interval, intersects $E$ in dimension $1/2$. Consider a cone of arbitrary small angular width $\theta$, centered around the $y-$axis, with vertex at the origin and pointing up.\footnote{Even though we originally restricted our set to the first quadrant, we can through a suitable rotation of coordinates ensure that the cone under consideration is symmetric about the y axis, and lies above the x axis.} In this case, for some large $ k_0 >0$, \footnote{We choose a scale $2^{2^{k}}$ instead of $2^k$ since if we have a range of $k$ values from 1 to some $k_0$, the levels are so sparse that only the last level $k_{0}$ is relevant when counting the points for the mass dimension till height $2^{2^{k_0}}$. With the scaling $2^{k_0}$ we would need to add up all the points in all the lower levels as well.} for each $k\geq k_0$, we fill the annular region inside the cone between the heights $2^{2^{k+1}}$ and $2^{2^{k}}+2^{2^{k+1}}$ with all the points belonging to $\mathbb{Z}^{2}$ within the annular region. In this way, we would have a set with mass dimension $3/2$ since we effectively have the area $ \approx (\theta)\cdot 4^{2^{k}}\cdot 2^{2^{k}}=\theta. 2^{3.2^{k}}$ covered just above the height $2^{2.2^{k}}$, by the integers, the mass dimension of the set being given by:
    
    \[ \lim_{k \to \infty} \frac{\log(\theta. 2^{3.2^{k}})}{\log(2^{2.2^{k}})}=\frac{3}{2} \]
    
    However, every tube contains at most  $2\cdot 2^{2^{k}}$ many points just above the height $2^{2^{k+1}}$ , and thus the mass dimension of each of these tubes is $1/2$.
    
    \item This is a second example of a set $E$ that has mass dimension $1$ and where an open cone of tubes has mass dimension $0$. We try to modify the previous example in order to reduce the dimension of $E$ from $3/2$ to $1$, however as it turns out in this specific case, each tube then only contains a finite number of points. Here assume that the cone is bounded by the angles $\pi/2$ and $(\pi/2- \theta)$ for some small $\theta$. Consider the set $E$ where starting with $k\geq0$ just above the height $2^{2^{k+1}}$, we have an annular integer grid chunk of height $2^{2^{k}}$, and radial width $\theta \cdot 2^{2^{k}}$ from the left edge of the cone, so that the area covered is exactly $\theta\cdot 2^{2^{k+1}}$ . This covers an angular width of $(\theta \cdot 2^{2^{k}})/2^{2^{k+1}}=\theta/2^{2^{k}}$. We do this for successive values of $k$ starting at $k=0$, at each step subtending an angle just to the right of the angle subtended in the previous case. At the height $2^{2^{k+1}}$ we put a square annular chunk as before of total area $2^{2^{k+1}}$, and in the process the total angle subtended as $k\to \infty$ is $\theta\cdot \sum\limits_{k=1}^{\infty} 1/2^{2^{k}}<\theta$. It is clear that this set $E$ has mass dimension $1$ since there are heights $2^{2^{k+1}}$ just above which an integer grid of area approximately $\theta \cdot 2^{2^{k+1}}$ is filled, but each tube with the right edge within the cone either has a finite positive number of points in it, or is empty, and thus has dimension 0. 
    
    \item   If we tried to avoid the infinite progression of angles that we got in the previous example, we can fix an angular width $\theta/2^{2^{k_0}}$ for each of the levels from heights $2^{2^{k_{0}+1}}$ to the height $2^{2^{k_{0}+ 2^{2^{k_{0}}}}}$, and put the annular chunk of integers of area $\theta\cdot 2^{2^{k_0 +1}}$ just above the height $2^{2^{k_{0}+1}}$ as before, so at this height the dimension of the set is close to $1$. 
    
    At the height $2^{2^{k_0 +2}}$, we put a chunk of width $\theta\cdot(1/2^{2^{k_{0}}})2^{2^{k_{0}+2}}=\theta\cdot 2^{3.2^{k_{0}}}$ and height $2^{2^{k_0 +1}}$ and thus an area of $\theta\cdot 2^{(2^2 +2-1).2^{k_{0}}}$. In the next stage ,we have a covered area of $\theta. 2^{7.2^{k_{0}}}.2^{2^{k_0 +2}}=\theta.2^{(2^3 +2^2 -1).2^{k_0}}$. just above the height $2^{2^{k_0 +3}}$. Continuing this way, in the end we have above the height $2^{2^{k_{0}+2^{2^{k_{0}}}}}$ , a covered area of  $\theta. 2^{2^{k_{0}}.(2^{2^{2^{k_0}}}+2^{2^{2^{k_0 }}-1}-1)}$. Thus the expression in the mass dimension (prior to taking the limit) when considering a box till this last height, would be:
    
    \[  \frac{\log(\theta\cdot 2^{2^{k_{0}}\cdot(2^{2^{2^{k_0}}}+2^{2^{2^{k_0 }}-1}-1)}}{\log(2^{2^{k_{0}+2^{2^{2^{k_{0}}}}}})}= \frac{\log(\theta)+\log(2). 2^{k_{0}}\cdot (2^{2^{2^{k_0}}}+(\frac{1}{2})2^{2^{2^{k_0 }}}-1)}{(\log 2\cdot 2^{k_{0}}\cdot 2^{ 2^{2^{k_{0}}}})}  \]
    
    When $k_{0}$ is large enough, this expression tends to $3/2$ from below and thus we have a set of dimension almost $3/2$ at this height. As we take the limit of $k_0 \to \infty$, the dimension of the set becomes exactly $3/2$. In this case, we don't have an infinite progression of angles; the full angle $\theta$ is covered after $2^{2^{k_{0}}}$ steps, however, the price to pay is that the dimension of $E$ increases to 3/2. Above the height $2^{2^{k_0 +2^{2^{k_{0}}}}}$ we can repeat the same process infinitely often, and then we end up having a set of dimension $3/2$ and where each of the tubes has dimension $1/2$ like before.
    

    
\end{enumerate}

\bigskip

In the Examples 3 and 5, we can clearly also modify the examples in a simple way so that instead of having sets $E$ of mass dimension $3/2$ and $E\cap t_{u,v}$ having mass dimension $1/2$, the sets and tubes have mass dimension $(1+\eta)$ and $\eta$ respectively, for any $0<\eta<1$.

\section{Proof of the slicing theorem.}

Here we finally prove \cref{thm:strongerslicing}. As noted in the introduction, \cref{thm:slicing}, the main slicing theorem, is an immediate consequence of this, upon integrating over the $v$ parameters.

\bigskip

\textit{Preliminaries:} We denote by $\mu$ the Lebesgue measure on the $u$ parameter space. We will always be working with a fixed value $\Tilde{v}$ on the $v$ parameter space, and will suppress the dependence of the measure $\mu$ on the fixed $\Tilde{v}$, which will be understood from the context.

We need to show there cannot exist for any fixed $\Tilde{v}$, any subset $U$ of the $u$ parameter  space with $\mu(U)>0$ so that

\[\overline{D}(E \cap t_{u,\Tilde{v}})> max(0, \overline{D}(E)-1 ), \  \text{when} \ (u,\Tilde{v})\in U.\]

Suppose to the contrary that there does exist some $\Tilde{v}$, and a subset $L$ of $u$ values with $\mu(L)>0$, so that for all $(u,\Tilde{v}),\ \text{with} \  u\in L$, we have $dim(E \cap t_{u,\Tilde{v}})> max(0, dim(E)-1 )$. We will show that this implies there exists a subset $E' \subset E$ with $dim(E')>dim(E)$ which is a contradiction.    

Given the above assumption, we must have some fixed $\psi>0$ so that there is some subset $L'\subset L$ with $\mu(L')>0$ where $L'$ is the set of exceptions where

\[\overline{D}(E\cap t_{u,\tilde{v}}) \geq max(0, \overline{D}(E)-1)+ \psi \ \text{for} \ u \in L'.\] 

Henceforth, we work with this set $L'$ and just call it $L$ itself. 

For any $\epsilon>0$, there is an open subset $O$ with $ L \subset O$ and $\mu(O-L)< \epsilon$.  Thus, less than an $\epsilon$ width of the angles within this open set are not covered by angles in $L$.\footnote{It's not necessary to make $\epsilon$ as small as possible; the cone is just a convenient set within which all the exceptional tubes lie. }

Without loss of generality, we can consider $O$ to be an open interval. This will correspond to an open cone pointing up and to the right and eventually intersecting the first quadrant where the $1$-separated set $E$ lies. It is enough to consider the case $\Tilde{v}=0$. The argument for any other $\Tilde{v}$ follows in a similar manner, the only difference being there is some finite initial region within the cone that does not intersect the first quadrant and which is thus empty; this initial empty space makes no difference as far as the mass dimension is concerned. 

For convenience, we rotate our coordinates so that the cone corresponding to $O$, containing this exceptional set of angles, points vertically upward and is symmetric about the $y$ axis. \footnote{This simplification helps later for calculating the mass dimension by taking boxes whose edges are parallel to the axes and hence also cutting perpendicularly across the cone, something that we are implicitly always doing throughout in our arguments.}

We denote the square box of length $n$ that is symmetric about the y axis and that lies in the upper half plane, as $B_{n}$, with one of the horizontal sides of $B_{n}$ lying along the $x$-axis. We denote the square box of length $n$, but whose one side is slanted parallel to the direction of the tube $t_{(u,0)}$ (i.e one side with slope $-(1/u)$), symmetric about the tube $t_{(u,0)}$, the other side of the box $B_{n}(u)$ lies along the line $y=ux$ symmetrically about the line $y=-\frac{1}{u}x$, as $B_{n}(u)$. This is illustrated in Figure 2.

\begin{figure}
\centering

\begin{tikzpicture}[x=1cm, y=1cm, semitransparent, scale=0.7]

\draw [<->] (0,-3) -- (0,6) node[above] {$y$} ;
\draw [<->] (-6,0) -- (6,0) node[right] {$x$} ;
\draw [-,,style=thick] (-2,2) -- (2,-2);
\draw [-,,style=thick] (2,-2) -- (6,2);
\draw [-, style=thick] (6,2)--(2,6) node[near end]{$B_{n}(u)$};
\draw [-, style=thick] (2,6)--(-2,2);
\draw[<->,style=thick] (2.3,-2.3)--(6.3,1.7)  node[below, pos=0.5] {$n$};
\draw[-,style=very thin] (0,0)--(5,5);
\draw[-, style=dashed] (-0.2,0.2)--(4.8,5.2);
\draw[-, style=dashed] (5,5)--(4.8,5.2) node[above] {$t_{(u,0)}$};
\draw[->, style= dashed] (0,0)--(7,4) node[right]{$l_{1}$}; 
\draw[->, style= dashed] (0,0)--(-7,4) node[left]{$l_{2}$};
\draw (-0.7,0.4) .. controls +(up:0.15cm) and +(right:0.5cm) .. (0.7,0.4)
\foreach \p in {0.5} {node[above,pos=\p]{$\theta$}};

\end{tikzpicture}
\caption{The tube $t_{(u,0)}$ and the `box' $B_{n}(u)$ are illustrated; here $u\in O $ and the lines $l_1$ and $l_2$ are the delimiting lines of the open cone of width $\theta$ corresponding to $O$.}
\end{figure}





Thus all the tubes under consideration are such that their right edges pass through the origin, and except for a small area $O(1)$ just near the origin, the tubes entirely lie within the cone.

\begin{lemma}\label{lemma1}

We write $\alpha=\text{max}(0,\text{dim}(E)-1)$. For any exceptional tube $t_{u,0}$ with $u\in L\cap O$, there are heights $m_{l}(u)$, with $l\in \N$, such that from the origin till $m_{l}(u)$, the upper half of the tube contains an exceptionally high number of points. Precisely,

\begin{equation}\label{eq:onlyeq}
    |(B_{m_{l}(u)}(u) \setminus B_{m_{l}(u)/2}(u))\cap E\cap t_{(u,0)}|>\Big(\frac{m_{l}(u)}{2}\Big)^{\alpha+ \psi/2}
\end{equation}
\end{lemma}

\begin{proof}[Proof of \cref{lemma1}:] For any exceptional tube $t_{u,0}$, start with any arbitrary level $n_{1}(u)$ large enough so that $n_{1}(u)^{\alpha +\psi }>K\cdot n_{1}(u)^{\alpha+\psi/2}$ and 

\[ \frac{\log|B_{n_{1}(u)}(u) \cap E \cap t_{(u,0)}|}{\log (n_{1}(u))}\geq \alpha +\psi \]  that is, 

\[ |B_{n_{1}(u)}(u)\cap E \cap t_{(u,0)}|\geq  n_{1}(u)^{\alpha + \psi} \]

where $K=\frac{1}{2^{\alpha}-1}$.

Choose the next level $n_{2}(u)$ so that $n_{2}(u) \geq e^{n_{1}(u)}$ and large enough that $n_{2}(u)^{\alpha +\psi }>K(u)\cdot n_{2}(u)^{\alpha+\psi/2}+\log n_{2}(u)$,  and where 

 \[ |B_{n_{2}(u)}(u) \cap E \cap t_{(u,0)}|\geq  n_{2}(u)^{\alpha + \psi} \] 

This implies there has to be some $k_{2}(u)\in \mathbb{N}$, so that $n_{2}(u)/2^{k_{2}(u)}> n_{1}(u) $ and so that,


\[ |\big(B_{n_{2}(u)/2^{k_{2}(u)-1}}(u) \setminus B_{n_{2}(u)/2^{k_{2}(u)}}(u)\big)  \cap E \cap t_{(u,0)}| > \Big(\frac{n_{2}(u)}{2^{k_{2}(u)}}\Big)^{\alpha+\psi/2} \]

 since otherwise, the total number of points within the tube from the origin to the $n_{2}$'th level is bounded above by:
 
\[ \sum\limits_{k=1}^{\infty} \Big(\frac{n_{2}(u)}{2^{k}}\Big)^{\alpha+\psi/2} + n_{1}(u)\leq K \cdot n_{2}(u)^{\alpha+\psi}+\log(n_{2}(u)), \] 

and we get a contradiction given our previous assumption.

Now we relabel the level $n_{2}(u)$ by the level $m'_{2}(u):=n_{2}(u)/2^{k_{2}(u)-1}$ we got from above. Next, by induction, we would find a level $n_{3}(u)$ and a positive integer $k_{3}(u)$ so that $n_{3}(u)/2^{k_{3}(u)}>n_{2}(u)$ and that 

\[ | \big(B_{n_{3}(u)/2^{k_{3}(u)-1}}(u) \setminus B_{n_{3}(u)/2^{k_{3}(u)}}(u)\big)  \cap E \cap t_{(u,0)}| > \Big(\frac{n_{3}(u)}{2^{k_{3}(u)}}\Big)^{\alpha+\psi/2}  .\]

Again we replace the level $n_{3}(u)$ by the level $m'_{3}(u):=n_{3}(u)/2^{k_{3}(u)-1}$. We further iterate the levels $n_{l}(u)$, for $l\geq 4$, by induction in the same way and find the corresponding levels $m'_{l}(u)$ so that 
\begin{equation}\label{eq:onlyeq}
    |(B_{m'_{l}(u)}(u) \setminus B_{m'_{l}(u)/2}(u))\cap E\cap t_{(u,0)}|>\Big(\frac{m'_{l}(u)}{2}\Big)^{\alpha+ \psi/2}
\end{equation}

 So for any exceptional tube $t_{u,0}$ within the cone with $u\in L$, we have found heights $m'_{l}(u)$, with $l\in \N$, such that from the origin till $m'_{l}(u)$, the upper half of the tube contains an exceptionally high number of points. For each tube $t_{(u,0)}$ with $u\in L$ and for each $l\in \N$, we inductively define the least integer greater than $m_{l-1}(u)$ that satisfies \cref{eq:onlyeq} (and which clearly exists) as $m_{l}(u)$.
 \end{proof}

\bigskip

Next, we have the following result:

\begin{lemma}\label{lemma2}
We have a subset of $L$ of positive measure on which for any fixed $l \in \N$, the set of values taken by the function $u \to m_l(u)$ is bounded.
\end{lemma}

We rename this aforementioned subset of $L$ as $L$ itself.

\begin{proof}[Proof of \cref{lemma2}:]
For a specific $u_0$, for any $l\in \mathbb{N}$, consider the level $m_{l}(u_0)$ that has been attained along the tube $t_{(u_0,0)}$. Since we have a $1$-separated set, till the height $n_{k}$ there are only finitely many points within the cone, which thus subtend finitely many angles at the origin. There are also a finite number of points contained in $t_{u_{0},0}$. We can thus slightly perturb $t_{(u_{0},0)}$ to the left by a certain finite amount, and no point will exit or enter the tube when this perturbation happens, when we restrict to looking at points only till the height $m_{l}(u_0)$. Once a point \big(till height $m_{l}(u_0)$ \big) exits the tube for the first time and we are at some perturbed angle $\tilde{u_{0}}\in L$ near $u_{0}$, while no new point has entered the tube, the $m_{l}(\tilde{u_{0}})$ value has to be equal or higher than $m_{l}(u_0)$: it cannot decrease since all the points in $t_{(\tilde{u_{0}},0)}$ till height $m_{l}(u)$ are also in $t_{(u_0,0)}$ and so $m_{l}(u_{0})$ value would have to be lower, which by definition is not possible. Similarly, it is verified that if a new point enters while no point has left the tube as the tube is rotated to the left, then the new $m_l$ value has to be equal or lower.

The same argument applies when we rotate the tube slightly to the right and we are considering the finite height $m_{l}$ as before; the only difference being that for certain angles of the tube, it might happen that a new point lies on the right edge or on the left edge of the tube which immediately enter and leave the tube respectively, as the tube is perturbed to the right. There being only finitely many points till any height $m_{l}$, the total number of such special angles of the tube when the above happens, is at most countable. \footnote{It only matters that the function $u \mapsto m_{l}(u)$ changes in discrete steps.}

Thus in summary, for a fixed $l\in \N$, values taken by the function $u \mapsto m_{l}(u)$ as $u$ varies in $L$, is discrete, hence the total number of values taken by this function $u \mapsto m_{l}(u)$ for fixed $l$, is countable as $u$ varies in $L$.



If the total number of values taken by the function $u \mapsto m_{l}(u)$ is finite for some particular $l\in \N$, then we keep all these tubes under consideration. If the total number of values taken by $u \mapsto m_{l}(u)$ for some fixed $l$ is countably infinite, we enumerate these values in increasing order of their $m_{l}$ values, and we can choose any arbitrarily small enough $ \epsilon<\mu(L)$ and discard a total $\epsilon/2^{l}$ width of tubes in $L$, and so the values taken by the function in the remaining set of values in $L$ is bounded.\footnote{Once the range of the function $u \mapsto m_{l_0}(u)$ is countably infinite, the range of all the functions $u \mapsto m_{l}(u)$ for $l>l_0$ is also countably infinite.} We do this for every single integer value of $l$, with the widths taken out being in a geometric progression, so that eventually we have taken out a total width at most $\sum\limits_{l=1}^{\infty} \frac{\eps}{2^{l}}=\eps$ from $L$ and still we are left with a positive measure set $L''\subset L$, and we simply rename this $L$. \footnote{Note that in Section 3 we have outlined examples where the levels $m_{l}$ diverge to infinity; and indeed since the functions $u \to n_{k}(u)$ are locally constant, the number of such diverging limiting angles are also countable. For the purposes of the slicing theorem, we are excising arbitrarily small neighborhoods around these diverging points and don't require a more detailed study of the structure of this set of diverging points.}

\begin{figure}
\centering

\begin{tikzpicture}[x=1cm, y=1cm, semitransparent, scale=0.7]

\draw [->] (0,0) -- (0,6) node[above] {$y$} ;
\draw [<->] (-4,0) -- (4,0) node[right] {$x$} ;
\draw [-] (-1,1.5)--(-0.95,1.5);
\draw [-] (-0.95,2)--(-0.9,2);
\draw [-] (-0.9,1.7)--(-0.8,1.7);
\draw [-] (-0.8,2)--(-0.7,2);
\draw [-] (-0.7,2.5)--(-0.6,2.5);
\draw [-] (-0.6,1.5)--(-0.5,1.5);
\draw [-] (-0.5,3.2)--(-0.4,3.2);
\draw [-] (-0.4,4.1)--(-0.3,4.1);
\draw [-] (-0.3,5)--(-0.2,5);
\draw [-] (0,3)--(0.5,3);
\draw [-] (0.5,3.5)--(1.2,3.5);
\draw [-] (1.2,4)--(2,4);
\draw [-] (2,4.5)--(3,4.5);
\draw[-,style=dashed] (0,0)--(-0.2,6) node[left]{$l_3$};
\draw[->, style= dashed] (0,0)--(4,6) node[right]{$l_{1}$}; 
\draw[->, style= dashed] (0,0)--(-4,6) node[left]{$l_{2}$};

\end{tikzpicture}
\caption{For a fixed $l\in \mathbb{N}$, an example of an $m_l$ `profile' is shown with the $u$-coordinates restricted to the subset $ L\cap O$, where $O$ corresponds to the cone illustrated above. The function $u \to m_{l}(u)$ for a fixed $l$ is discrete as explained in the proof, thus takes countable many values. There may be instances of the values diverging, for example to some line $l_3$ from the left as shown above, in which case by removing from the domain of the function arbitrarily small neighborhoods around these diverging points (By removing a set of tubes  whose $u$ values lie in a set of measure $\eps/2^{l}$, the above function is bounded).}
\end{figure}

\bigskip

So finally we have a subset $L$ of positive measure on which all the $m_{l}(u)$ `profiles' are bounded.
\end{proof}

With \cref{lemma1} and \cref{lemma2} in hand, we complete the proof of \cref{thm:strongerslicing}.

\begin{proof}[Proof of \cref{thm:strongerslicing}]

We construct a subset $E'$ of $E$ which has dimension greater than the dimension of $E$ itself, which gives us the contradiction we need. As stated before, our cone points vertically upward and is centered around the y axis, with vertex at the origin.

Let $\mu(L)=\beta>0 $. For any $l\in \mathbb{N}$, let the supremum of the values of $m_{l}(u)$ as $u\in L$, be $e^{H_{l}}$. We start with a value $l_{0}$ so that the minimum height $h_{l_0}$ of the $m_{l_{0}}$ profile is some arbitrary large number, say 100, and from now on we consider $l\geq l_{0}$.\footnote{We can always find such a level $l_0$.}  We divide the range $[0,e^{H_{l}}]$ into $H_{l}$ distinct parts of equal height, so for $1 \leq i \leq \lceil H_l \rceil $, the $i$'th interval corresponds to the height range $[\frac{i-1}{H_l}e^{H_l},\frac{i}{H_l}e^{H_{l}}]$. By the pigeonhole principle, for any $l\in \N$, the pre-image of at least one of these $\lceil H_l \rceil$ height-ranges under the map $u \mapsto m_{l}(u), u\in L$, has measure greater than $\beta/\lceil H_l \rceil$. Call this preimage set $T_{l}\subset L$.  First we outline our argument when for any $l\in \N$ the height range of $T_{l}$ is some $[\frac{j-1}{H_l}e^{H_l},\frac{j}{H_l}e^{H_{l}}]$ for $2\leq j \leq \lceil H_l \rceil$. The case $j=1$ is treated later. 


Consider the arc $A_l$ inside the cone at height $\frac{j-1}{2H_l}e^{H_l}$ from the origin, and the subset $A'_l$ of $A_l$ which corresponds to the intersections of $A_l$ with tubes in $T_{l}$.\footnote{We mean intersections of the right edge of the tube with $A_{l}$, since the tubes are prametrized by the angles of the right edges.}. The set $A'_{l}$ has Lebesgue measure greater than $\frac{(j-1)\beta}{2( H_l +1)^{2}}e^{H_l}$, since the pre-image set has measure greater than $\beta/\lceil H_{l} \rceil$ . We split $A'_{l}$ into disjoint consecutive arcs each of Lebesgue measure 1, and color these disjoint lengths alternately blue and red. We consider the union of the blue arcs, and the union of the red arcs. Then the intersection of the tubes in $T_l$ with at least one of these two arcs has measure greater than or equal to $\frac{(j-1)\beta}{4(H_l +1)^{2}}e^{H_l}$. Call this the good arc. Thus, per unit length of the good arc if we just counted one exceptional tube, we have $\frac{(j-1)\beta}{4(H_l +1)^{2}}e^{H_l}$ many distinct tubes belonging to $T_{l}$. Call this set of tubes $T'_{l}$. These tubes have been chosen so that beyond the height $\frac{j-1}{2H_{l}}e^{H_{l}}$, they are all mutually disjoint, so when we add all the points within all these tubes, there is no double counting.

Because of the way the levels $m_{l}$ have been defined, we must have more than $\Big(\frac{(j-1)}{2(H_l +1)}e^{H_l}\Big)^{\alpha+\psi/2}$ many points in each of the tubes of $T'_{l}$, in the height range $[\frac{j-1}{2H_{l}}e^{H_{l}},\frac{j}{H_{l}}e^{H_{l}}]$. Thus in this height range, there are strictly more than  :
\begin{equation}\label{eq:equation1}
\frac{(j-1)\beta}{4(H_l +1)^{2}}e^{H_l}\cdot \Big(\frac{(j-1)}{2(H_l +1)}e^{H_l}\Big)^{\alpha + \psi/2} = e^{H_{l}\cdot (1+\alpha+ \psi/2)}\frac{\beta(j-1)^{1+\alpha+\psi/2}}{2^{2+\alpha+\psi/2}\cdot (H_l +1)^{2+ \alpha +\psi/2}}
\end{equation}

many points, where now: $1+\alpha+\psi/2 =\text{max}(1, \text{dim}(E))+\psi/2$.  Thus within this set, at height $\frac{j}{H_{l}}e^{H_{l}}$, the expression in the formula for mass dimension is:

\begin{equation}\label{eq:equation2}
\frac{\Bigg(\splitfrac{H_{l}(1+\alpha +\psi/2)+\log(\beta)+(1+\alpha +\psi/2)\log(j-1)}{-(2+\alpha+\psi/2)\log(2)-(2+\alpha+\psi/2)\log(H_{l}+1)}\Bigg)}{H_{l}+\log(j)-\log(H_{l})} 
\end{equation}

The leading terms in both the numerator and denominator above contain a factor of $H_{l}$, and when $H_{l}$ is sufficiently large, this expression approaches $(1+\alpha+\psi/2)$ which is strictly greater than $\text{dim}(E)$.


However if $j=1$, the above process fails since there is no lower bound to the $m_{l}$ profile, and we cannot choose a lower cut off height to perform the integration as we did above. In this case, we subdivide again and write $\frac{1}{H_{l}}e^{H_{l}}=e^{H^{(1)}_{l}}$ where $H^{(1)}_{l}=H_{l}-log(H_{l})< H_{l}$. Now, if for any $2\leq j \leq \lceil H_{l} \rceil$, we take out a pre-image set of measure at least $\beta/\lceil H_{l} \rceil^{3}$, it is clear from \cref{eq:equation2} that upon taking logarithms we still get a dimension greater than $\text{dim}(E)$ when the heights are sufficiently large, since the only change is an additive term in the coefficient of $\log(H_{l}+1)$ in the numerator.  


So we assume that for $2\leq j\leq \lceil H_{l} \rceil$, at most $\beta/\lceil H_{l} \rceil^{3} $ measure of the set has been taken out, and so we are left with at least $(1-\frac{(\lceil H_{l} \rceil -1)}{\lceil H_{l} \rceil^{3}})\beta >(1- \frac{1}{\lceil H_{l} \rceil^{2}})\beta $ measure for the pre-image set of the range $[0,\frac{1}{H_{l}}e^{H_{l}}]=[0, e^{H^{(1)}_{l}}]$.

Now we divide this range again as before, and then if there is some $j\geq 2$ for which the preimage set of $[\frac{j-1}{H^{(1)}_{l}}e^{H_{l}^{(1)}},\frac{j}{H^{(1)}_{l}}e^{H_{l}^{(1)}}]$ has measure at least $\beta/\lceil (H_{l}^{(1)}) \rceil^{3} $, then we can count the number of points as in the previous case. If not, we have a preimage set of measure at least $(1- \frac{1}{\lceil H_{l} \rceil^{2} }-\frac{1}{\lceil H_{l}^{1} \rceil^{2}})\beta$ for the range $[0, \frac{1}{H_{l}^{(1)}}e^{H_{l}^{(1)}}]=[0, e^{H_{l}^{(2)}}]$
where $H_{l}^{(2)}=H_{l}^{(1)}-\log(H_{l}^{(1)})<H_{l}^{(1)}$. Eventually we would reach the minimum height $h_{l}$  of this profile and still be left with a set of measure at least $(1-(\frac{1}{100^{2}}+....\frac{1}{H_{l}^{2}}))\beta>C\beta$ for some constant $C\approx (1-\frac{1}{100})$ independent of $l$.

At this very end for any $l\in \N$, we reach a final height range $H^{t(l)}_{l}$ for some $t(l)\in \mathbb{N}$, so that the minimum of the $m_{l}$ profile is some height $p(l)$ so that $p(l)> \frac{1}{H_{l}^{t(l)}}e^{H_{l}^{t(l)}}$ in which case $p(l)\in [\frac{j-1}{H_{l}^{t(l)}}e^{H_{l}^{t(l)}},\frac{j}{H_{l}^{t(l)}}e^{H_{l}^{t(l)}}]$, for some $j\geq 2$(If we were to hit the $j=1$ range here, then by definition we would reach a lower height range than $H^{t(1)}_{1}$). 

Thus, the total number of points in this height range is at least 

\begin{equation}\label{eq:equation3}
    \frac{(j-1)C\beta e^{H_{l}^{t(l)}}}{4(H_{l}^{t(l)}+1)}\cdot \Bigg(\frac{(j-1)e^{H_{l}^{t_{l}}}}{2(H_{l}^{t(l)}+1)} \Bigg)^{\alpha+\psi/2}
\end{equation}

From the construction, it is apparent that the function $l \mapsto H_{l}^{t(l)}$ is strictly increasing. Looking at the structure of \cref{eq:equation2}, we see that the analogous expression here would also only involve changes by some additive term in the coefficient of the $\log(H_{l}^{t_{l}})$ terms, and so as $H_{l}^{t(l)}\to \infty$, the mass dimension still approaches $1+\alpha +\psi/2$.


It is clear that in the worst case, even if we were always reduced to this lowest height range for each $l$, since the minimum height profiles $h_{l}\to \infty$ as $l\to \infty$, we always get an increasing sequence of squares $B_{n_{l}}$ with $n_{l}\to \infty$, and by considering its intersections with the cone ensure that $E$ has an exceptionally high dimension. For every $l$ we get the required height profile at which the mass dimension is close to $(\text{max}(1,\text{dim}(E))+\psi/2)$ and thus we can choose a subsequence from these profiles for each $l$ value, and conclude that the dimension is exactly $(\text{max}(1,\text{dim}(E))+\psi/2)$.\footnote{We will eventually always find a sub-sequence along which the dimension becomes exceptionally large.} Thus we constructed the set $E'\subset E$ whose mass dimension is greater than the mass dimension of $E$ itself, which is the required contradiction. 
\end{proof}

Now the proof of \cref{corollary:finalcor} follows immediately from \cref{thm:strongerslicing}:

\begin{proof}[Proof of \cref{corollary:finalcor}]

In this case, we take $E=A\times B$, a Cartesian product. We note that the set $\{(x,y)\in \R^{2}| y=\lfloor ux+v \rfloor\}$ for $u>0$ always lies inside a tube of width less than one. In fact, the vertical cross section of the set $\lfloor l_{u,v} \rfloor$ is always 1, and the width is clearly less than one. Thus the result follows from \cref{thm:slicing}. 
\end{proof}

\section{Future directions:}
In a companion manuscript, we intend to study the Hausdorff dimensions of the exceptional sets of $u$-parameters that violate the inequality of \cref{thm:strongerslicing}, for specific $1$-separated sets $E$. The proof of the \cref{thm:strongerslicing} does not follow along the lines of the proof of the classical Marstrand slicing theorem, as Examples 1 and 2 in Section 3 demonstrate, and the exceptional sets also do not behave in classically expected ways. 

One should also naturally be able to extend the slicing theorem to higher dimensions.

\section{Acknowedgements.} The author is grateful to Daniel Glasscock for suggesting this problem, and for helpful discussions. The author is thankful to Dmitry Kleinbock for general discussions, and to Hershdeep Singh for a discussion on Example 3 of Section 3.

\end{document}